\documentclass[11pt]{article}
\usepackage[numbers]{natbib}

\usepackage{amsmath,amssymb}
\usepackage{enumerate,hyperref}
\numberwithin{equation}{section}
\newtheorem{theorem}{Theorem}[section]
\newtheorem{lemma}{Lemma}[section]

\newtheorem{remark}{Remark}[section]
\newtheorem{proposition}{Proposition}[section]

\def\QED{\mbox{\rule[0pt]{1.5ex}{1.5ex}}}
\def\endproof{\hspace*{0.98\textwidth}~\QED\par\endtrivlist\unskip}

\usepackage{amsmath}

\newcommand{\cov}{\mathrm{Cov}}
\newcommand{\iid}{{i.i.d.}}

\newcommand{\RR}{\mathbb{R}}
\newcommand{\LL}{\mathbb{L}}

\title{\textbf{Adaptive estimation of an additive regression function from weakly dependent data}}
\author{Christophe Chesneau\footnote{LMNO, CNRS-Universit\'e de Caen, Campus II, Science 3, 14032, Caen, France.}, Jalal Fadili\footnote{GREYC, CNRS-ENSICAEN-Universi\'e de Caen, 6, Bd du Mar\'echal Juin, 14050 Caen Cedex, France.} and Bertrand Maillot$^*$}
\date{}

\begin{document}

\maketitle

\begin{abstract}
A $d$-dimensional nonparametric additive regression model with dependent observations is considered. Using the marginal
integration technique and wavelets methodology, we develop a new adaptive estimator for a component of the additive regression function. Its asymptotic properties are investigated via the minimax approach under the $\mathbb{L}_2$ risk over Besov balls. We prove that it attains a sharp rate of convergence which turns to be the one obtained in the $\iid$ case for the standard univariate regression estimation problem. 
\end{abstract}

\vspace{9pt}
\par
\vspace{9pt}
\par

\section{Introduction}
\label{sec:intro}
\subsection{Problem statement}
Let $d$ be a positive integer, $(Y_i,{\bf X}_i)_{i\in \mathbb{Z}}$
be a $\RR\times [0,1]^d$-valued stationary process on a
probability space $(\Omega,\mathcal{A},\mathbb{P})$ and $\rho$ be a
given real measurable function. The unknown regression function
associated to $(Y_i,{\bf X}_i)_{i\in \mathbb{Z}}$ and $\rho$ is
defined by  
\[
g({\bf x})=\mathbb{E}(\rho(Y)| {\bf X}={\bf x}),
\qquad{\bf x}=(x_1,\ldots,x_d)\in [0,1]^d.
\]
In the additive regression model, the function $g$ is considered to have an additive structure, i.e. there exist $d$ unknown real measurable functions
$g_1,\ldots,g_d$ and an unknown real number $\mu$ such that
\begin{eqnarray} \label{gg}
g({\bf x})= \mu+\sum_{\ell=1}^dg_{\ell}(x_{\ell}).
\end{eqnarray}
For any $\ell\in \{1,\ldots,d\}$, our goal is to estimate $g_{\ell}$
from $n$ observations $(Y_1,{\bf X}_1),\ldots,(Y_n,{\bf X}_n)$ of
$(Y_i,{\bf X}_i)_{i\in \mathbb{Z}}$.

\subsection{Overview of previous work}
When $(Y_i,{\bf X}_i)_{i\in \mathbb{Z}}$ is a $\iid$ 
process, this additive regression model becomes the standard
one. In such a case,  Stone in a series of papers \cite{stone1,stone2,stone3} proved that
$g$ can be estimated with the same rate of estimation error as in
the one-dimensional case. The estimation of the component $g_{\ell}$
has been investigated in several papers via various methods (kernel,
splines, wavelets, etc.). See e.g. \cite{buja}, \cite{hastie},
\cite{linton}, \cite{ops1,ops2}, \cite{am1}, \cite{am2}, \cite{sper},
\cite{zhang}, \cite{sar} and \cite{fanjiang}.

In some applications, the $\iid$ assumption on the
observations is too stringent. For this reason, some authors
have explored the estimation of $g_{\ell}$ in the dependent case.
When $(Y_i,{\bf X}_i)_{i\in \mathbb{Z}}$ is a strongly mixing
process, this problem has been addressed  by \cite {camlong},
\cite{deb0}, and results for continuous time processes under a strong
mixing condition have been obtained by \cite{deb1,deb2}. In
particular, they have developed non-adaptive kernel estimators for
$g_{\ell}$ and studied its asymptotic properties.

\subsection{Contributions}
To the best of our knowledge, adaptive estimation of $g_{\ell}$ for
dependent processes has been addressed only by \cite{gao}. The lack
of results for adaptive estimation in this context motivates this work.
To reach our goal, as in \cite{zhang}, we combine the marginal
integration technique introduced by \cite{newey} with wavelet
methods. We capitalize on wavelets to construct an adaptive thresholding estimator 
and show that it attains sharp rates of convergence under mild assumptions on the
smoothness of the unknown function. By adaptive, it is meant that the parameters of the estimator do not depend on the parameter(s) of the dependent process nor on those of the smoothness class of the function. In particular, this leads to a simple and easily implementable estimator.

More precisely, our wavelet estimator is based on term-by-term hard thresholding. 
The idea of this estimator is simple:
(i) we estimate the unknown wavelet coefficients of $g_{\ell}$  based on the
observations; (ii) then we select the greatest ones and ignore the others; (iii) and finally we
reconstruct the function estimate from the chosen wavelet coefficients on the
considered wavelet basis. Adopting the minimax point of view under
the $\mathbb{L}_2$ risk, we prove that our adaptive estimator attains a sharp rate of convergence over Besov balls which capture a variety of smoothness features in a function including spatially inhomogeneous behavior. The attained rate corresponds to the optimal one in the $\iid$ case for the univariate regression
estimation problem (up to an extra logarithmic term).

\subsection{Paper organization}
The rest of the paper is organized as follows. Section~\ref{sec:assump} presents our assumptions on the model. In Section~\ref{sec:wav}, we describe wavelet bases on $[0,1]$, Besov balls and tensor product wavelet bases on $[0,1]^d$. Our wavelet hard thresholding estimator is detailed in Section~\ref{sec:estim}.
Its rate of convergence under the $\mathbb{L}_2$ risk over Besov balls is established in Section~\ref{sec:minimax}. Section~\ref{sec:relprior} provides a discussion of the relation of our result with respect to prior work. The proofs are detailed in Section~\ref{sec:proofs}.

\section{Notations and assumptions}
\label{sec:assump}
In this work, we assume the following on our model:

\paragraph{Assumptions on the variables.}
\begin{itemize}
\item For any $i\in \{1,\ldots,n\}$, we set ${\bf X}_i=(X_{1,i},\ldots,X_{d,i})$. We suppose that
\begin{itemize}
\item for any $i\in \{1,\ldots,n\}$, $X_{1,i},\ldots,X_{d,i}$ are identically distributed with the common distribution $\mathcal{U}([0,1])$,
\item ${\bf X}_1,\ldots,{\bf X}_n$ are  identically distributed with the common known density $f$.
\end{itemize}
\item We suppose that the following identifiability condition is satisfied: for any $\ell\in \{1,\ldots,d\}$ and $i\in \{1,\ldots,n\}$, we have
\begin{eqnarray}\label{iden}
\mathbb{E}(g_{\ell}(X_{\ell,i}))=0.
\end{eqnarray}
\end{itemize}

\paragraph{Strongly mixing assumption.}
Throughout this work, we use the strong mixing dependence structure on $(Y_i,{\bf X}_i)_{i\in \mathbb{Z}}$. For any $m\in \mathbb{Z}$, we define the $m$-th strongly mixing coefficient of $(Y_i,{\bf X}_i)_{i\in \mathbb{Z}}$ by
\begin{eqnarray}\label{stronglyy}
\alpha_m=\sup_{(A,B)\in \mathcal{F}^{(Y,{\bf X})}_{-\infty,0}\times \mathcal{F}^{(Y,{\bf X})}_{m,\infty}}\left| \mathbb{P}(A\cap B)-\mathbb{P}(A)\mathbb{P}(B)\right|,
\end{eqnarray}
where  $\mathcal{F}^{(Y,{\bf X})}_{-\infty,0}$ is the $\sigma$-algebra generated by $\ldots, (Y_{-1},{\bf X}_{-1}),(Y_0,{\bf X}_0)$ and $\mathcal{F}^{(Y,{\bf X})}_{m,\infty}$ is the $\sigma$-algebra generated by $(Y_m,{\bf X}_m), (Y_{m+1},{\bf X}_{m+1}),\ldots$ ~.

We suppose that there exist two constants $\gamma>0$ and $c>0$ such that, for any integer $m\ge 1$,
\begin{eqnarray}\label{dedans}
\alpha_m\le \gamma \exp(-c m).
\end{eqnarray}
Further details on
strongly mixing dependence can be found in \cite{bradl},
\cite{withers}, \cite{doukhan}, \cite{moda} and \cite{car}.

\paragraph{\bf Boundedness assumptions.}
\begin{itemize}
\item We suppose that $\rho \in \LL_1(\RR) \cap \LL_\infty(\RR)$, i.e. there exist constants $C_1>0$ and $C_2$ (supposed known) such that
\begin{eqnarray}
\label{refff2}
& \int_{-\infty}^{\infty}|\rho(y)|dy\le C_1, \\
\label{refff}
\text{and} & \sup_{y\in \RR}|\rho(y)|\le C_2.
\end{eqnarray}
\item We suppose that there exists a known constant $c>0$ such that
\begin{eqnarray}\label{fal}
\inf_{{\bf x}\in [0,1]^d}f({\bf x})\ge c.
 \end{eqnarray}
\item For any $m\in \{1,\ldots,n\}$, let $f_{(Y_0,{\bf X}_{0},Y_m,{\bf X}_{m})}$ be the density of $(Y_0,{\bf X}_{0},Y_m,{\bf X}_{m})$, $f_{(Y_0,{\bf X}_{0})}$ the density of $(Y_0,{\bf X}_{0})$ and, for any  $(y,{\bf x},y_*,{\bf x}_*)\in \RR\times [0,1]^d\times \RR\times [0,1]^d$,
\begin{eqnarray}\label{chemin}
\lefteqn{h_m(y,{\bf x},y_*,{\bf x}_*) =}& & \nonumber \\
&
& f_{(Y_0,{\bf X}_{0},Y_m,{\bf X}_{m})}(y,{\bf x},y_*,{\bf x}_*)-f_{(Y_0,{\bf X}_{0})}(y,{\bf x})f_{(Y_0,{\bf X}_{0})}(y_*,{\bf x}_*). \nonumber \\
&
&
\end{eqnarray}
We suppose that there exists a known constant $C>0$ such that
\begin{eqnarray}\label{boundw}
\sup_{m\in\{1,\ldots,n\}}\sup_{(y,{\bf x},y_*,{\bf x}_*)\in \RR\times [0,1]^d\times \RR\times [0,1]^d}|h_m(y,{\bf x},y_*,{\bf x}_*) |\le C.
\end{eqnarray}
\end{itemize}
Such boundedness assumptions are standard for the estimation of $g_{\ell}$ from a strongly mixing process. See e.g. \cite{deb1,deb2}.

\section{Wavelets and Besov balls}
\label{sec:wav}
\subsection{Wavelet bases on $[0,1]$}
Let $R$ be a positive integer. We consider an orthonormal wavelet basis generated by dilations and translations of the scaling and wavelet functions $\phi$ and $\psi$ from the Daubechies family $\mathrm{db}_{2R}$. In particular, $\phi$ and $\psi$ have compact supports and unit $\LL_2$-norm, and $\psi$ has $R$ vanishing moments, i.e. for any $r\in \{0,\ldots,R-1\}$, $\int x^r\psi(x)dx=0$.

Define the scaled and translated version of $\phi$ and $\psi$
\[
\phi_{j,k}(x)=2^{j/2}\phi(2^jx-k),  \qquad \psi_{j,k}(x)=2^{j/2}\psi(2^jx-k).
\]
Then, with an appropriate treatment at the boundaries, there exists an integer $\tau$ satisfying $2^\tau\ge 2R$ such that, for any integer $j_*\ge \tau$, the collection
\[
\{ \phi_{j_*,k}(.), \ k\in \{0,\ldots,2^{j_*}-1\}; \ \psi_{j,k}(.); \ j \in \mathbb{N}-\{0,\ldots, j_*-1\} ,\ k\in \{0,\ldots,2^{j}-1\}\},
\] 
is an orthonormal basis of $\mathbb{L}_2(\lbrack 0,1 \rbrack)=\{h : [0,1]\rightarrow \RR; \ \ \int_{0}^{1}h^2(x)dx <\infty\}$. See \cite{cohen2,mallat}.

Consequently, for any integer $j_*\ge \tau$, any $h\in \mathbb{L}_2(\lbrack 0,1 \rbrack)$ can be expanded into a
wavelet series as
$$h(x)= \sum_{k=0}^{2^{j_*}-1}\alpha_{j_*,k}\phi_{j_*,k}(x)  +\sum_{j= j_*}^{\infty}  \sum_{ k=0}^{2^j-1}\beta_{j,k}\psi_{j,k}(x),\qquad x\in [0,1],$$
where
\begin{eqnarray}\label{avion1}
\alpha_{j,k}=\int_{0}^{1}h(x)\phi_{j,k}(x)dx, \qquad \beta_{j,k}=\int_{0}^{1}h(x)\psi_{j,k}(x)dx.
\end{eqnarray}

\subsection{Besov balls} 
As is traditional in the wavelet estimation literature, we will investigate the performance of our
estimator by assuming that the unknown function to be estimated belongs to a Besov ball. The Besov norm for a
function can be related to a sequence space norm on its wavelet coefficients. 
More precisely, let $M>0$, $s \in (0,R)$, $ p \ge 1$ and $q \ge 1$.
A function $h$ in $\mathbb{L}_2(\lbrack 0,1 \rbrack)$ belongs to $ \mathbf{B}^s_{p,q}(M)$ if, and only if, there exists a constant $M^*>0$ (depending on $M$) such that the associated wavelet coefficients \eqref{avion1} satisfy
\begin{eqnarray*}
 \left(\sum_{j=\tau}^{\infty} \left(2^{j(s+1/2-1/p )}\left(\sum_{k=0}^{2^{j}-1}|\beta_{j,k}|^{p}\right)^{1/p}\right)^{q}\right)^{1/q} \le  M^*.
\end{eqnarray*}
In this expression, $s$ is a smoothness parameter and $p$ and $q$ are norm parameters. Besov spaces include many
traditional smoothness spaces. For particular choices of $s$, $ p$ and $q$, Besov balls contain the standard H\"older and Sobolev balls. See \cite{meyer}.

\subsection{Wavelet tensor product bases on $[0,1]^d$}
For the purpose of this paper, we will use compactly supported tensor product wavelet bases on $[0,1]^d$ based on the Daubechies family. Let us briefly recall their construction. For any ${\bf x}=(x_1, \ldots,x_d)\in [0,1]^d$, we construct a scaling function 
\[
\Phi({\bf x})=\prod_{v=1}^d \phi (x_v) ~,
\]
and $2^d-1$ wavelet functions
\begin{eqnarray*}
\Psi_u({\bf x})= \left\{
\begin{aligned}
& \psi(x_{u})\prod_{\underset{v\not = u}{v=1}}^d\phi(x_{v})& & {\text{when $u\in \{1,\ldots,d\}$}} ,\\
& \prod_{v\in A_u}\psi(x_{v})\prod_{v\not \in A_u}\phi(x_{v}) &  & {\text{when $u\in \{d+1,\ldots,2^d-1\}$,}}
\end{aligned}
\right.
\end{eqnarray*}
where $(A_u)_{u\in \{d+1,\ldots,2^d-1\}}$ forms the set of all non void subsets of $\{1,\ldots,d\}$ of cardinality greater or equal to $2$.

For any integer $j$ and any ${\bf k}=(k_1,\ldots,k_d)$, define the translated and dilated versions of $\Phi$ and $\Psi_u$ as
\begin{align*}
\Phi_{j,{\bf k}}({\bf x})&=2^{jd/2}\Phi(2^jx_1-k_1, \ldots,2^jx_d-k_d), \\
\Psi_{j,{\bf k},u}({\bf x})&=2^{jd/2}\Psi_{u}(2^jx_1-k_1, \ldots,2^jx_d-k_d), ~ \text{for any $u\in \{1,\ldots,2^d-1\}$}.
\end{align*}

Let $D_j=\{0,\ldots,2^j-1\}^d$. Then, with an appropriate treatment at the boundaries, there exists an integer $\tau$ such that the collection 
\[
\{\Phi_{\tau,{\bf k}},
{\bf k} \in D_{\tau}; \ (\Psi_{j,{\bf k},u})_{u\in \{1,\ldots,2^d-1\}}, \ \ \ j\in \mathbb{N}-\{0,\ldots,\tau-1\}, \
{\bf k}\in  D_{j}\}
\] 
forms an orthonormal basis of $\mathbb{L}_2(\lbrack 0,1 \rbrack^d)=\{h : [0,1]^d\rightarrow \RR; \ \ \int_{[0,1]^d}h^2({\bf x})d{\bf x} <\infty\}$.

For any integer $j_*$ such that $j_* \ge \tau$, a function $h \in \mathbb{L}_2(\lbrack 0,1 \rbrack^d)$ can be expanded into a
wavelet series as
$$
h({\bf x})= \sum_{{\bf k}\in  D_{j_*}}\alpha_{j_*,{\bf k}}\Phi_{j_*,{\bf k}}({\bf x})  +\sum_{u=1}^{2^d-1}\sum_{j= j_*}^{\infty}  \sum_{{\bf k}\in  D_j}\beta_{j,{\bf k},u}\Psi_{j,{\bf k},u}({\bf x}), \qquad {\bf x}\in [0,1]^d,
$$
where
\begin{eqnarray}\label{avion}
\alpha_{j,{\bf k}}=\int_{[0,1]^d}h({\bf x})\Phi_{j,{\bf k}}({\bf x})d{\bf x}, \qquad \beta_{j,{\bf k},u}=\int_{[0,1]^d}h({\bf x})\Psi_{j,{\bf k},u}({\bf x})d{\bf x}.
\end{eqnarray}

\section{The estimator}
\label{sec:estim}
\subsection{Wavelet coefficients estimator}
The following proposition provides a wavelet decomposition of $g_{\ell}$ based on the ``marginal integration'' method (introduced by \cite{newey}) and a tensor product wavelet basis on $[0,1]^d$.

\begin{proposition}\label{woo} Suppose that \eqref{iden} holds. Then, for any $j_*\ge \tau$ and $\ell\in \{1,\ldots,d\}$, we can write
\[
g_{\ell}(x)= \sum_{k=1}^{2^{j_*}-1}a_{j_*,k,\ell}\phi_{j_*,k}(x)  +\sum_{j= j_*}^{\infty}  \sum_{k=1}^{2^j-1}b_{j,k,\ell}\psi_{j,k}(x) - \mu,\qquad x\in [0,1],
\]
where
\begin{eqnarray}\label{coef1}
a_{j,k,\ell} = a_{j,k_{\ell},\ell}=2^{-j(d-1)/2}\int_{[0,1]^d}g({\bf x})\sum_{{\bf k}_{-\ell}\in D_{j}^*}\Phi_{j,{\bf k}}({\bf x})d{\bf x}, \\
\label{coef2}
 b_{j,k,\ell} = b_{j,k_{\ell},\ell} = 2^{-j(d-1)/2}\int_{[0,1]^d}g({\bf x})\sum_{{\bf k}_{-\ell}\in D_{j}^*}\Psi_{j,{\bf k},\ell}({\bf x})d{\bf x},
\end{eqnarray}
and ${\bf k}_{-\ell}=(k_1,\ldots,k_{\ell-1},k_{\ell+1},\ldots,k_d)$ and $D_j^*=\{0,\ldots,2^j-1\}^{d-1}$.
\end{proposition}

\begin{remark}
Due to the definitions of $g$ and properties of $\Psi_{j,{\bf k},\ell}$, $b_{j,k,\ell}$ is nothing but the wavelet coefficient of $g_{\ell}$, i.e.
\begin{eqnarray}\label{berti}
 b_{j,k,\ell}=\int_{0}^{1}g_{\ell}(x)\psi_{j,k}(x)dx=\beta_{j,k}.
\end{eqnarray}
\end{remark}

Proposition \ref{woo} suggests that a first step to estimate $g_{\ell}$ should consist in estimating the unknown coefficients $a_{j,k,\ell}$ \eqref{coef1} and $b_{j,k,\ell}$ \eqref{coef2}. To this end, we propose the following coefficients estimators
\begin{eqnarray}\label{ru}
\widehat a_{j,k,\ell}=\widehat a_{j,k_{\ell},\ell}=2^{-j(d-1)/2}\frac{1}{n}\sum_{i=1}^n\frac{\rho(Y_i)}{f({\bf X}_i)}\sum_{{\bf k}_{-\ell}\in D_{j}^*}\Phi_{j,{\bf k}}({\bf X}_i)
\end{eqnarray}
and
\begin{eqnarray}\label{quier}
\widehat b_{j,k,\ell}=\widehat b_{j,k_{\ell},\ell}=2^{-j(d-1)/2}\frac{1}{n}\sum_{i=1}^n\frac{\rho(Y_i)}{f({\bf X}_i)}\sum_{{\bf k}_{-\ell}\in D_{j}^*}\Psi_{j,{\bf k},\ell}({\bf X}_i).
\end{eqnarray}

These estimators enjoy powerful statistical properties. Some of them are collected in the following propositions.

\begin{proposition}[Unbiasedness]\label{unbiased} Suppose that \eqref{iden} holds. For any $j\ge \tau$, $\ell\in \{1,\ldots,d\}$ and $k\in \{0,\ldots,2^j-1\}$, $\widehat a_{j,k,\ell}$ and $\widehat b_{j,k,\ell}$ in \eqref{ru} and \eqref{quier} are unbiased estimators of $a_{j,k,\ell}$ and $b_{j,k,\ell}$ respectively.
\end{proposition}

\begin{proposition}[Moment inequality I]\label{var22} Suppose that the assumptions of Section~\ref{sec:assump} hold. Let $j\ge \tau$ such that $2^j\le n$, $k\in \{0,\ldots,2^j-1\}$, $\ell\in \{1,\ldots,d\}$. 
Then there exists a constant $C>0$ such that
\begin{eqnarray*}
\mathbb{E}\left( (\widehat a_{j,k,\ell}-a_{j,k,\ell} )^2\right)\le   C \frac{1}{n}, \qquad  \mathbb{E}\left( (\widehat b_{j,k,\ell}-b_{j,k,\ell} )^2\right)\le   C \frac{1}{n}.
\end{eqnarray*}
\end{proposition}

\begin{remark}
In the proof of Proposition \ref{var22}, we only need to have the existence of two constants $C>0$ and  $q\in (0,1)$ such that  $\sum_{m=1}^{n}m^{q}\alpha_{m}^{q}\le C<\infty$. This latter inequality is obviously satisfied by \eqref{dedans}.
\end{remark}

\begin{proposition}[Moment inequality II]\label{var2} 
Under the same assumptions of Proposition~\ref{var22}, there exists a constant $C>0$ such that
\begin{eqnarray*}
\mathbb{E}\left( (\widehat b_{j,k,\ell}-b_{j,k,\ell} )^4\right)\le   C\frac{2^j}{n}.
\end{eqnarray*}
\end{proposition}

\begin{proposition}[Concentration inequality]\label{var3} Suppose that the assumptions of Section~\ref{sec:assump} hold. Let $j\ge \tau$ such that $2^j\le n/(\ln n)^3$, $k\in \{0,\ldots,2^j-1\}$, $\ell\in \{1,\ldots,d\}$ and $\lambda_n=(\ln n/n)^{1/2}$. Then there exist two constants $C>0$ and $\kappa>0$ such that
\begin{eqnarray*}
 \mathbb{P}\left( |\widehat b_{j,k,\ell}-b_{j,k,\ell} | \ge \kappa \lambda_n/2\right)\le     C\frac{1}{n^4}.
\end{eqnarray*}
\end{proposition}

\subsection{Hard thresholding estimator}
We now turn to the estimator of $g_\ell$ from $\widehat a_{j,k,\ell}$ and $\widehat b_{j,k,\ell}$ as introduced in \eqref{ru} and \eqref{quier}. Towards this goal, we will only keep the significant wavelet coefficients that are above a certain threshold according to the hard thresholding rule, and then reconstruct from these coefficients. In a compact form, this reads
\begin{eqnarray} \label{hinz2}
\widehat g_{\ell}(x)= \sum_{k=0}^{2^{\tau}-1}\widehat a_{\tau,k,\ell} \phi_{\tau,k}(x)+\sum_{j=\tau}^{j_1}\sum_{k=0}^{2^j-1}\widehat b_{j,k,\ell}{\bf 1}_{\left \lbrace |\widehat b_{j,k,\ell}|\ge \kappa\lambda_n \right \rbrace }\psi_{j,k}(x)-\widehat \mu,
\end{eqnarray}
where
\begin{eqnarray}\label{ahm}
\widehat \mu =\frac{1}{n}\sum_{i=1}^n\rho(Y_i).
\end{eqnarray}
In \label{hinz2}, $j_1$ is the resolution level satisfying $ 2^{j_1}=[n/(\ln n)^3]$,
$\kappa$ is a large enough constant (the one in Proposition \ref{var3}) and
\begin{eqnarray*}
\lambda_n= \sqrt{\frac{\ln n }{n}}.
\end{eqnarray*}

Note that,  due to the assumptions on the model, our wavelet hard
thresholding estimator \eqref{hinz2} is simpler than the one of
\cite{zhang}. 

\section{Minimax upper-bound result}
\label{sec:minimax}
Theorem~\ref{theo2} below investigates the minimax rates of convergence attained by $\widehat g_{\ell}$ over Besov balls
under the $\LL_2$ risk.
\begin{theorem}\label{theo2} 
Let $\ell \in \{1,\ldots,d\}$. Suppose that the assumptions of Section~\ref{sec:assump} hold. Let $\widehat g_{\ell}$ be the estimator given in \eqref{hinz2}. Suppose that $g_{\ell}\in \mathbf{B}^s_{p,q}(M)$ with $q\ge 1$, \{$p \ge 2$ and $s \in (0,R)$\} or \{$p \in [1,2)$ and $s \in (1/p,R)$\}. Then there exists a constant $C>0$ such that
\[
\mathbb{E}\left(\int_{0}^{1}(\widehat g_{\ell}(x)-g_{\ell}(x))^2dx\right)\le C \left( \frac{\ln n}{n}\right)^{2s/(2s+1)}.\]
\end{theorem}
The proof of Theorem \ref{theo2} is based on a suitable decomposition of the $\mathbb{L}_2$ risk and the statistical properties of \eqref{ru} and \eqref{quier} summarized in Propositions~\ref{unbiased}, \ref{var22}, \ref{var2} and \ref{var3} above.

\section{Relation to prior work}
\label{sec:relprior}
The rate $(\ln n/n)^{2s/(2s+1)}$ is, up to an extra logarithmic term, known to be the optimal one for the standard one-dimensional regression model with uniform random design. See e.g. \cite{hardle} and \cite{tsybakov}.

Theorem~\ref{theo2} provides an ``adaptive contribution'' to the results of \cite {camlong}, \cite{deb0} and  \cite{deb1,deb2}. Furthermore, if we confine ourselves to the $\iid$ case, we recover a similar result to \cite[Theorem 3]{zhang} but without the condition $s>\max (d/2,d/p)$. The price to pay is more restrictive assumptions on the model ($\rho$ is bounded from above, the density of ${\bf X}$ is known, etc.). Additionally, our estimator has a more straightforward and friendly implementation than the one in \cite{zhang}.

\section{Proofs}
\label{sec:proofs}
In this section, the quantity $C$ denotes any constant that does not depend on $j$, $k$ and $n$. Its value may change from one term to another and may depends on $\phi$ or $\psi$.

\subsection{Technical results on wavelets}
\paragraph{Proof of Proposition \ref{woo}.}
Because of \eqref{refff}, we have $g\in \mathbb{L}_2([0,1]^d)$. For any $j_*\ge \tau$, we can expand $g$ on our wavelet-tensor product basis as

\begin{eqnarray}\label{bur1}
g({\bf x})= \sum_{{\bf k}\in  D_{j_*}}\alpha_{j_*,{\bf k}}\Phi_{j_*,{\bf k}}({\bf x})  +\sum_{u=1}^{2^d-1}\sum_{j= j_*}^{\infty}  \sum_{{\bf k}\in  D_j}\beta_{j,{\bf k},u}\Psi_{j,{\bf k},u}({\bf x}), \qquad {\bf x}\in [0,1]^d
\end{eqnarray}

where $$\alpha_{j,{\bf k}}=\int_{[0,1]^d}g({\bf x})\Phi_{j,{\bf k}}({\bf x})d{\bf x}, \qquad\beta_{j,{\bf k},u}=\int_{[0,1]^d}g({\bf x})\Psi_{j,{\bf k},u}({\bf x})d{\bf x}.$$

Moreover, using the ``marginal integration'' method based on \eqref{iden}, we can write
\begin{eqnarray}\label{bur2}
g_{\ell}(x_{\ell})=\int_{[0,1]^{d-1}}g({\bf x}) \prod_{\underset{v\not = \ell}{v=1}}^d dx_{v}- \mu, \qquad x_{\ell}\in [0,1].
\end{eqnarray}
Since $\int_{0}^{1}\phi_{j,k}(x)dx=2^{-j/2}$ and $\int_{0}^{1}\psi_{j,k}(x)dx=0$, observe that
$$\int_{[0,1]^{d-1}}\Phi_{j_*,{\bf k}}({\bf x}) \prod_{\underset{v\not = \ell}{v=1}}^d dx_{v}=2^{-j_*(d-1)/2}\phi_{j_*,k_{\ell}}(x_{\ell})$$ and
\begin{eqnarray*}
\int_{[0,1]^{d-1}}\Psi_{j,{\bf k},u}({\bf x})\prod_{\underset{v\not = \ell}{v=1}}^d dx_{v}= \left\{
\begin{aligned}
& 2^{-j(d-1)/2}\psi_{j,k_{\ell}}(x_{\ell})& & {\text{if $u=\ell$}} ,\\
& 0&  & {\text{otherwise.}}
\end{aligned}
\right.
\end{eqnarray*}
Therefore, putting \eqref{bur1} in \eqref{bur2} and writing $x=x_{\ell}$, we obtain
\begin{eqnarray*}
g_{\ell}(x) = \sum_{{\bf k}\in  D_{j_*}}2^{-j_*(d-1)/2}\alpha_{j_*,{\bf k}}\phi_{j_*,k_{\ell}}(x)  +\sum_{j= j_*}^{\infty}  \sum_{{\bf k}\in  D_j}2^{-j(d-1)/2}\beta_{j,{\bf k},\ell}\psi_{j,k_{\ell}}(x) - \mu.
\end{eqnarray*}
Or, equivalently,
$$g_{\ell}(x)= \sum_{k=1}^{2^{j_*}-1}a_{j_*,k,\ell}\phi_{j_*,k}(x)  +\sum_{j= j_*}^{\infty}  \sum_{k=1}^{2^j-1}b_{j,k,\ell}\psi_{j,k}(x) - \mu,$$
where
\begin{eqnarray*}
a_{j,k,\ell}=a_{j,k_{\ell},\ell}=2^{-j(d-1)/2}\int_{[0,1]^d}g({\bf x})\sum_{{\bf k}_{-\ell}\in D_{j}^*}\Phi_{j,{\bf k}}({\bf x})d{\bf x}
\end{eqnarray*}
and
\begin{eqnarray*}
 b_{j,k,\ell}= b_{j,k_{\ell},\ell}=2^{-j(d-1)/2}\int_{[0,1]^d}g({\bf x})\sum_{{\bf k}_{-\ell}\in D_{j}^*}\Psi_{j,{\bf k},\ell}({\bf x})d{\bf x}.
\end{eqnarray*}
Proposition \ref{woo} is proved.

\endproof

\begin{proposition}\label{mout}
For any $\ell\in \{1,\ldots,d\}$, $j\ge \tau$ and $k=k_{\ell}\in \{0,\ldots,2^j-1\}$, set
$$h^{(1)}_{j,k}({\bf x})=\sum_{{\bf k}_{-\ell}\in D_{j}^*}\Phi_{j,{\bf k}}({\bf x}), \ \ \ \ h^{(2)}_{j,k}({\bf x})=\sum_{{\bf k}_{-\ell}\in D_{j}^*}\Psi_{j,{\bf k},\ell}({\bf x}), \qquad {\bf x}\in  [0,1]^d.$$
Then there exists a constant $C>0$ such that, for any $a\in \{1,2\}$,
 $$\sup_{{\bf x}\in [0,1]^d}|h^{(a)}_{j,k}({\bf x})|\le C2^{jd/2}, \qquad \int_{ [0,1]^d}|h^{(a)}_{j,k}({\bf x})|d{\bf x}\le C2^{-j/2} 2^{j(d-1)/2}$$ and $$\int_{ [0,1]^d} (h^{(a)}_{j,k}({\bf x}))^2d{\bf x}=2^{j(d-1)}.$$
\end{proposition}

\begin{proof}
\begin{itemize}
\item Since $\sup_{x\in [0,1]}|\phi_{j,k}(x)|\le C2^{j/2}$ and $\sup_{x\in [0,1]}\sum_{k=0}^{2^j-1}|\phi_{j,k}(x)|\le C2^{j/2}$, we obtain
\begin{eqnarray*}
\sup_{{\bf x}\in [0,1]^d}|h^{(1)}_{j,k}({\bf x})| = (\sup_{x\in [0,1]}|\phi_{j,k}(x)|)  \left( \sup_{x\in [0,1]}\sum_{k=0}^{2^j-1}|\phi_{j,k}(x)|\right)^{d-1} \le  C2^{jd/2}.
\end{eqnarray*}
\item Using $\int_{0}^{1}|\phi_{j,k}(x)|dx=C2^{-j/2}$, we obtain
\begin{eqnarray*}
\int_{ [0,1]^d}|h^{(1)}_{j,k}({\bf x})|d{\bf x}& \le &   \left(\int_{0}^{1}|\phi_{j,k}(x)|dx\right) \left(\sum_{k=0}^{2^j-1}\int_{0}^{1}|\phi_{j,k}(x)|dx\right)^{d-1} \\
& =
& C 2^{-j/2} 2^{j(d-1)/2}.
\end{eqnarray*}
\item Since, for any $(u_{\bf k})_{{\bf k}\in D_j}$, $\int_{ [0,1]^d}\left(\sum_{{\bf k}\in D_{j}}u_{\bf k}\Phi_{j,{\bf k}}({\bf x})\right)^2 d{\bf x}=\sum_{{\bf k}\in D_{j}}u_{\bf k}^2$, we obtain
\begin{eqnarray*}
\int_{ [0,1]^d}(h^{(1)}_{j,k}({\bf x}))^2d{\bf x} = \int_{ [0,1]^d}\left(\sum_{{\bf k}_{-\ell}\in D_{j}^*}\Phi_{j,{\bf k}}({\bf x})\right)^2 d{\bf x}= 2^{j(d-1)}.
\end{eqnarray*}
Proceeding in a similar fashion, using $\sup_{x\in [0,1]}|\psi_{j,k}(x)|\le C2^{j/2}$, $\int_{0}^{1}|\psi_{j,k}(x)|dx=C2^{-j/2}$ and, for any $(u_{\bf k})_{{\bf k}\in D_j}$, $\int_{ [0,1]^d}\left(\sum_{{\bf k}\in D_{j}}u_{\bf k}\Psi_{j,{\bf k},\ell}({\bf x})\right)^2 d{\bf x}=\sum_{{\bf k}\in D_{j}}u_{\bf k}^2$, we obtain the same results for $h^{(2)}_{j,k}$.
\end{itemize}
This ends the proof of Proposition \ref{mout}.

\end{proof}

\subsection{Statistical properties of the coefficients estimators}
\paragraph{Proof of Proposition \ref{unbiased}.} 
We have
\begin{eqnarray*}
\mathbb{E}(\widehat a_{j,k,\ell})& =& 2^{-j(d-1)/2}\mathbb{E}\left(\frac{\rho(Y_1)}{f({\bf X}_1)}\sum_{{\bf k}_{-\ell}\in D_{j}^*}\Phi_{j,{\bf k}}({\bf X}_1)\right)\\
& =& 2^{-j(d-1)/2}\mathbb{E}\left(\mathbb{E}(\rho(Y_1)| {\bf X}_1)\frac{1}{f({\bf X}_1)}\sum_{{\bf k}_{-\ell}\in D_{j}^*}\Phi_{j,{\bf k}}({\bf X}_1)\right)\\
& =& 2^{-j(d-1)/2}\mathbb{E}\left(\frac{g({\bf X}_1)}{f({\bf X}_1)}\sum_{{\bf k}_{-\ell}\in D_{j}^*}\Phi_{j,{\bf k}}({\bf X}_1)\right)\\
& =
& 2^{-j(d-1)/2}\int_{[0,1]^d}\frac{g({\bf x})}{f({\bf x})}\sum_{{\bf k}_{-\ell}\in D_{j}^*}\Phi_{j,{\bf k}}({\bf x}) f({\bf x})d{\bf x}\\
& =
& 2^{-j(d-1)/2}\int_{[0,1]^d}g({\bf x})\sum_{{\bf k}_{-\ell}\in D_{j}^*}\Phi_{j,{\bf k}}({\bf x}) d{\bf x}= a_{j,k,\ell}.
\end{eqnarray*}
Proceeding in a similar fashion, we prove that $\mathbb{E}(\widehat b_{j,k,\ell})=b_{j,k,\ell}$.


\endproof

\paragraph{Proof of Proposition \ref{var22}.} 
For the sake of simplicity, for any $i\in \{1,\ldots,n\}$, set
$$Z_i=\frac{\rho(Y_i)}{f({\bf X}_i)}\sum_{{\bf k}_{-\ell}\in D_{j}^*}\Phi_{j,{\bf k}}({\bf X}_i).$$
Thanks to Proposition \ref{unbiased}, we have
\begin{eqnarray} \label{mom}
\mathbb{E}\left( (\widehat a_{j,k,\ell}-a_{j,k,\ell} )^2\right) = \mathbb{V}(\widehat a_{j,k,\ell}) = 2^{-j(d-1)}\frac{1}{n^2} \mathbb{V}\left(\sum_{i=1}^n Z_i\right).
\end{eqnarray}
An elementary covariance decomposition gives
\begin{eqnarray}\label{bouche}
 \mathbb{V}\left( \sum_{i=1}^n Z_i\right)& =&   n\mathbb{V}\left( Z_1\right)+2\sum_{v=2}^n\sum_{u=1}^{v-1}\cov \left(   Z_v,Z_u\right)\nonumber\\
& \le
& n\mathbb{V}\left( Z_1\right) + 2 \left|\sum_{v=2}^n\sum_{u=1}^{v-1}\cov \left( Z_v,Z_u\right)\right|.
\end{eqnarray}
Using \eqref{refff}, \eqref{fal} and Proposition \ref{mout}, we have
\begin{eqnarray}\label{bouche2}
\mathbb{V}\left( Z_1\right) & \le &  \mathbb{E}(Z_1^2)\le \frac{\sup_{y\in \RR}\rho^2(y)}{\inf _{{\bf x}\in [0,1]^d} f({\bf x})} \mathbb{E}\left( \frac{1}{f({\bf X}_1)}\left(\sum_{{\bf k}_{-\ell}\in D_{j}^*}\Phi_{j,{\bf k}}({\bf X}_1)\right)^2\right) \nonumber\\
& \le
&C \int_{[0,1]^d}\frac{1}{f({\bf x})}\left(\sum_{{\bf k}_{-\ell}\in D_{j}^*}\Phi_{j,{\bf k}}({\bf x})\right)^2f({\bf x})d{\bf x}\nonumber \\
& =
& C \int_{[0,1]^d}\left(\sum_{{\bf k}_{-\ell}\in D_{j}^*}\Phi_{j,{\bf k}}({\bf x})\right)^2d{\bf x}= C2^{j(d-1)}.
\end{eqnarray}
It follows from the stationarity of $(Y_i,{\bf X}_i)_{i\in \mathbb{Z}}$ and $2^j\le n$ that
\begin{eqnarray}\label{t}
\left|\sum_{v=2}^n\sum_{u=1}^{v-1}\cov \left(  Z_v,Z_u\right)\right| = \left|\sum_{m=1}^{n}(n-m)\cov \left(  Z_0,Z_m\right)\right|\le  R_1+R_2,
\end{eqnarray}
where
\[
R_1=n\sum_{m=1}^{2^j-1}\left|\cov \left( Z_0,Z_m\right)\right|, \ \ \ \ \ R_2=n\sum_{m=2^j}^{n}\left|\cov \left( Z_0,Z_m\right)\right|.
\]
It remains to bound $R_1$ and $R_2$.

\begin{enumerate}[(i)]
\item {\bf Bound for $R_1$.} Let, for any $(y,{\bf x},y_*,{\bf x}_*)\in \RR\times [0,1]^d\times \RR\times [0,1]^d$, $h_m(y,{\bf x},y_*,{\bf x}_*)$ be \eqref{chemin}. Using \eqref{boundw},  \eqref{refff2} and Proposition \ref{mout}, we obtain
\begin{eqnarray*}
\lefteqn{\left|\cov \left(Z_0,Z_m\right)\right|} & & \\
& =
& \bigg|\int_{-\infty}^{\infty}\int_{[0,1]^d}\int_{-\infty}^{\infty}\int_{[0,1]^d}h_m(y,{\bf x},y_*,{\bf x}_*)\times \\
&
&  \left(\frac{\rho(y)}{f({\bf x})}\sum_{{\bf k}_{-\ell}\in D_{j}^*} \Phi_{j,{\bf k}}({\bf x})\frac{\rho(y_*)}{f({\bf x}_*)}\sum_{{\bf k}_{-\ell}\in D_{j}^*} \Phi_{j,{\bf k}}({\bf x}_*)\right) dyd{\bf x}dy_*d{\bf x}_*\bigg|\\
& \le
&\int_{-\infty}^{\infty}\int_{[0,1]^d}\int_{-\infty}^{\infty}\int_{[0,1]^d} |h_m(y,{\bf x},y_*,{\bf x}_*)| \times\\
&
& \left|\frac{\rho(y)}{f({\bf x})}\right| \left|\sum_{{\bf k}_{-\ell}\in D_{j}^*} \Phi_{j,{\bf k}}({\bf x}) \right|  \left|\frac{\rho(y_*)}{f({\bf x}_*)}\right| \left|\sum_{{\bf k}_{-\ell}\in D_{j}^*} \Phi_{j,{\bf k}}({\bf x}_*)\right| dyd{\bf x}dy_*d{\bf x}_*\\
& \le
& C \left(  \int_{-\infty}^{\infty} |\rho(y)|dy\right)^2 \left(  \int_{[0,1]^d} \left|\sum_{{\bf k}_{-\ell}\in D_{j}^*}\Phi_{j,{\bf k}}({\bf x})\right|d{\bf x}\right)^2\le C2^{-j} 2^{j(d-1)}.
\end{eqnarray*}
Therefore
\begin{eqnarray}\label{t1}
R_1\le Cn2^{-j} 2^{j(d-1)}2^{j}=Cn2^{j(d-1)}.
\end{eqnarray}

\item {\bf Bound for $R_2$.}
By the Davydov inequality for strongly mixing processes (see \cite{davy}), for any $q\in (0,1)$, we have
\begin{eqnarray*}
\lefteqn{\left|\cov \left(   Z_0,Z_m\right)\right|\le   10 \alpha_{m}^q \left( \mathbb{E}\left( |Z_0|^{2/(1-q)}\right)\right)^{1-q}} & &  \\
& \le
& 10 \alpha_{m}^q \left(\frac{\sup_{y\in \RR}|\rho(y)|}{\inf _{{\bf x}\in [0,1]^d} f({\bf x})}\sup_{{\bf x}\in [0,1]^d}\left|\sum_{{\bf k}_{-\ell}\in D_{j}^*}\Phi_{j,{\bf k}}({\bf x})\right| \right)^{2q} \left( \mathbb{E}( Z_0^{2})\right)^{1-q}.
\end{eqnarray*}
By \eqref{refff}, \eqref{fal} and Proposition \ref{mout}, we have
\begin{eqnarray*}
\frac{\sup_{y\in \RR}|\rho(y)|}{\inf _{{\bf x}\in [0,1]^d} f({\bf x})}  \sup_{{\bf x}\in [0,1]^d}\left|\sum_{{\bf k}_{-\ell}\in D_{j}^*}\Phi_{j,{\bf k}}({\bf x})\right| & \le & C\sup_{{\bf x}\in [0,1]^d}\left|\sum_{{\bf k}_{-\ell}\in D_{j}^*}\Phi_{j,{\bf k}}({\bf x})\right| \\
& \le
& C 2^{jd/2}.
\end{eqnarray*}
By \eqref{bouche2}, we have
$$\mathbb{E}\left( Z_0^{2}\right)\le C2^{j(d-1)}.$$
Therefore
\begin{eqnarray*}
\left|\cov \left( Z_0,Z_m\right)\right|\le  C2^{q j}2^{j(d-1)}\alpha_{m}^q.
\end{eqnarray*}
Observe that $\sum_{m=1}^{\infty}m^q\alpha_{m}^q= \gamma^q \sum_{m=1}^{\infty}m^q exp(-c qm)<\infty$. Hence
\begin{eqnarray}\label{t2}
R_2\le Cn2^{q j}2^{j(d-1)}\sum_{m=2^j}^{n}\alpha_{m}^q\le Cn2^{j(d-1)}\sum_{m=2^j}^{n}m^q\alpha_{m}^q\le Cn2^{j(d-1)}.
\end{eqnarray}
\end{enumerate}

Putting \eqref{t}, \eqref{t1} and \eqref{t2} together, we have
\begin{eqnarray}\label{yeux}
\left|\sum_{v=2}^n\sum_{u=1}^{v-1}\cov \left(   Z_v,Z_u\right)\right| \le   C n2^{j(d-1)}.
\end{eqnarray}
Combining \eqref{mom}, \eqref{bouche}, \eqref{bouche2} and \eqref{yeux}, we obtain
\begin{eqnarray*}
\mathbb{E}\left( (\widehat a_{j,k,\ell}-a_{j,k,\ell} )^2\right)\le C2^{-j(d-1)}\frac{1}{n^2}n 2^{j(d-1)}= C \frac{1}{n}.
\end{eqnarray*}
Proceeding in a similar fashion, we prove that
\begin{eqnarray*}
\mathbb{E}\left( (\widehat b_{j,k,\ell}-b_{j,k,\ell} )^2\right)\le  C \frac{1}{n}.
\end{eqnarray*}
This ends the proof of Proposition \ref{var22}.

\endproof

\paragraph{Proof of Proposition \ref{var2}.} 
It follows from \eqref{refff}, \eqref{fal} and Proposition \ref{mout} that
\begin{eqnarray*}
|\widehat b_{j,k,\ell}| & \le &2^{-j(d-1)/2}\frac{1}{n}\sum_{i=1}^n \frac{|\rho(Y_i) |}{|f({\bf X}_i)|}\left|\sum_{{\bf k}_{-\ell}\in D_{j}^*}\Psi_{j,{\bf k},\ell}({\bf X}_i)\right|\\
& \le
&  2^{-j(d-1)/2}\frac{\sup_{y\in \RR}|\rho(y)|}{\inf _{{\bf x}\in [0,1]^d} f({\bf x})}  \sup_{{\bf x}\in [0,1]^d}\left|\sum_{{\bf k}_{-\ell}\in D_{j}^*}\Psi_{j,{\bf k},\ell}({\bf x})\right| \\
& \le
& C 2^{-j(d-1)/2} 2^{jd/2}=C2^{j/2}.
\end{eqnarray*}
Because of \eqref{refff}, we have $\sup_{{\bf x}\in [0,1]^d}|g({\bf x})|\le C$. It follows from Proposition \ref{mout} that
\begin{eqnarray}\label{plus}
|b_{j,k,\ell}| & \le & 2^{-j(d-1)/2}\int_{[0,1]^d}|g({\bf x})| \left|\sum_{{\bf k}_{-\ell}\in D_{j}^*}\Psi_{j,{\bf k},\ell}({\bf x})\right|d{\bf x}\nonumber\\
& \le
& C2^{-j(d-1)/2}\int_{[0,1]^d} \left|\sum_{{\bf k}_{-\ell}\in D_{j}^*}\Psi_{j,{\bf k},\ell}({\bf x})\right|d{\bf x}\nonumber\\
& \le
& C 2^{-j(d-1)/2}2^{-j}2^{jd/2}=C2^{-j/2}.
\end{eqnarray}
Hence
\begin{eqnarray}\label{bu}
|\widehat b_{j,k,\ell}-b_{j,k,\ell}| \le |\widehat b_{j,k,\ell}|+|b_{j,k,\ell}| \le C 2^{j/2}.
\end{eqnarray}
It follows from \eqref{bu} and Proposition \ref{var22}  that
\begin{eqnarray*}
\mathbb{E}\left( (\widehat b_{j,k,\ell}-b_{j,k,\ell})^4\right) \le  C  2^{j}\mathbb{E}\left( (\widehat b_{j,k,\ell}-b_{j,k,\ell})^2\right)\le C\frac{2^j}{n}.
\end{eqnarray*}

The proof of Proposition \ref{var2} is complete.

\endproof

\paragraph{Proof of Proposition \ref{var3}.} 
Let us first state a Bernstein inequality for exponentially strongly mixing process. 
\begin{lemma}[\cite{lieb2}]\label{bnh}
Let $\gamma>0$, $c>0$ and $(Y_i)_{i\in \mathbb{Z}}$ be a
stationary process with the $m$-th strongly mixing coefficient $\alpha_m$ \eqref{stronglyy}. Let $n$ be a positive
integer, $h : \RR \rightarrow \mathbb{C}$ be a measurable
function and, for any $i\in \mathbb{Z}$, $U_i=h(Y_i)$. We assume
that  $\mathbb{E}(U_1)=0$ and there exists a constant $M>0$
satisfying $|U_1| \le M$. Then, for any $m\in \{1,\ldots,[n/2]\}$ and
$\lambda>0$, we have
\[
\mathbb{P}\left( \left | \frac{1}{n}\sum_{i=1}^n U_i \right | \ge \lambda \right) \le 4\exp\left(-\frac{\lambda^2  n}{ 16(D_m/m+\lambda M m/{3})} \right)+32\frac{M}{\lambda}n\alpha_m,
\]
where $D_m=\max_{l\in \{1,\ldots,2m\}}\mathbb{V} \left( \sum_{i=1}^l U_i\right)$. 
\end{lemma}

We now apply this lemma by setting for any $i\in \{1,\ldots,n\}$, 
\[
U_{i}=2^{-j(d-1)/2}\frac{\rho(Y_i)}{f({\bf X}_i)}\sum_{{\bf k}_{-\ell}\in D_{j}^*}\Psi_{j,{\bf k},\ell}({\bf X}_i)-b_{j,k,\ell}.
\]
Then we can write
$$\widehat b_{j,k,\ell}-b_{j,k,\ell}=\frac{1}{n}\sum_{i=1}^n U_i.$$
So
\begin{eqnarray*}
\mathbb{P}\left( |\widehat b_{j,k,\ell}-b_{j,k,\ell}| \ge \kappa\lambda_n/2 \right)=\mathbb{P}\left( \left|\frac{1}{n}\sum_{i=1}^n U_i\right| \ge \kappa\lambda_n/2 \right),
\end{eqnarray*}
where $U_{1},\ldots,U_{n}$ are identically distributed, depend on $(Y_i,{\bf X}_i)_{i\in \mathbb{Z}}$ satisfying \eqref{dedans},
\begin{itemize}
\item by Proposition \ref{unbiased}, we have $\mathbb{E}(U_1)=0$,
\item using arguments similar  to Proposition \ref{var22} with $l$ instead of $n$, we prove that 
$$\mathbb{V} \left( \sum_{i=1}^l U_i\right)\le C l.$$
Hence $D_m=\max_{l\in \{1,\ldots,2m\}}\mathbb{V} \left( \sum_{i=1}^l U_i\right)\le  Cm$. 
\item proceeding in a similar fashion to \eqref{bu}, we obtain $|U_1| \le C2^{j/2}$.
\end{itemize}
Lemma \ref{bnh} applied with the random variables $U_1,\ldots,U_{n}$, $\lambda=\kappa \lambda_n/2$, $\lambda_n=(\ln n/n)^{1/2}$, $m=u\ln n$ with $u>0$ (chosen later),  $M=C2^{j/2}$, $2^j\le n/ (\ln n)^3$ and \eqref{dedans} gives
\begin{eqnarray*}
\lefteqn{\mathbb{P}\left( |\widehat b_{j,k,\ell}-b_{j,k,\ell}| \ge \kappa\lambda_n/2 \right)} & & \\
& \le
&  C\left(\exp\left(-C\frac{ \kappa^2\lambda_n^2 n }{1+\kappa \lambda_n mM} \right)+\frac{M}{\lambda_n}n\exp(-c m)\right)  \nonumber\\
& \le
& C\left(\exp\left(-C\frac{ \kappa^2\ln n}{1+{\kappa  u2^{j/2}\ln n(\ln n/n)^{1/2}}} \right)+\frac{2^{j/2}}{(\ln n/n)^{1/2}}n\exp(-c u\ln n) \right) \nonumber\\
& \le
& C\left( n^{-C \kappa^2/(1+\kappa u)} +n^{1-cu}\right).
\end{eqnarray*}
Therefore, for large enough $\kappa$ and $u$, we have
\begin{eqnarray*}
\mathbb{P}\left( |\widehat b_{j,k,\ell}-b_{j,k,\ell}| \ge \kappa\lambda_n/2 \right)\le C \frac{1}{n^4}.
\end{eqnarray*}
This ends  the proof of Proposition \ref{var3}.

\endproof

\subsection{Proof of Theorem~\ref{theo2}}

Using Proposition \ref{woo}, we have
\begin{eqnarray*}
\lefteqn{\widehat g_{\ell}(x)-g_{\ell}(x)} & & \\
 & =
& \sum_{k=0}^{2^{\tau}-1}(\widehat \alpha_{\tau,k,\ell}-\alpha_{\tau,k,\ell})\phi_{\tau,k}(x)+\sum_{j=\tau}^{j_1}\sum_{k=0}^{2^{j}-1}(\widehat b_{j,k,\ell}{\bf 1}_{\left\lbrace  |\widehat b_{j,k,\ell}|\ge \kappa\lambda_n \right\rbrace}-b_{j,k,\ell})\psi_{j,k}(x)  \\
& -
& \sum_{j=j_1+1}^{\infty}\sum_{k=0}^{2^{j}-1}b_{j,k,\ell}\psi_{j,k}(x)-(\widehat \mu-\mu).
\end{eqnarray*}
Using the elementary inequality: $(x+y)^2\le 2(x^2+y^2)$, $(x,y)\in\RR^2$, and the orthonormality of the wavelet basis, we have
\begin{eqnarray}\label{rou}
\mathbb{E}\left(\int_{0}^{1}(\widehat g_{\ell}(x)-g_{\ell}(x))^2dx\right)\le 2(T+U+V+W),
\end{eqnarray}
where
\begin{eqnarray*}
T=\mathbb{E}( (\widehat \mu-\mu)^2), &\qquad& U=\sum_{k=0}^{2^{\tau}-1}\mathbb{E}\left((\widehat \alpha_{\tau,k,\ell}-\alpha_{\tau,k,\ell})^2\right), \\
V=\sum_{j=\tau}^{j_1}\sum_{k=0}^{2^{j}-1}\mathbb{E}\left((\widehat b_{j,k,\ell}{\bf 1}_{\left\lbrace  |\widehat b_{j,k,\ell}|\ge \kappa\lambda_n \right\rbrace}-b_{j,k,\ell})^2\right), &\qquad& W=\sum_{j=j_1+1}^{\infty}\sum_{k=0}^{2^{j}-1}b_{j,k,\ell}^2.
\end{eqnarray*}

\begin{enumerate}[(i)]
\item {\bf Bound for $T$.} We proceed as in the proof of Proposition \ref{var22}. By \eqref{iden}, we have $\mathbb{E}(\rho(Y_1))=\mu$. Thanks to the stationarity of $(Y_i)_{i\in\mathbb{Z}}$, we have
\begin{eqnarray*}
T = \mathbb{V}(\widehat \mu)\le  \frac{1}{n}\mathbb{V}(\rho(Y_1))+2\frac{1}{n}\sum_{m=1}^{n}\left|\cov \left(  \rho(Y_0),\rho(Y_m)\right)\right|.
\end{eqnarray*}
Using \eqref{refff}, the Davydov inequality (see \cite{davy}) and \eqref{dedans}, we obtain
\begin{eqnarray}\label{del}
T\le C \frac{1}{n} \left(1+\sum_{m=1}^{n}\alpha_{m}^q\right)\le C\frac{1}{n}\le C\left( \frac{\ln n}{n}\right)^{2s/(2s+1)} .
\end{eqnarray}

\item \noindent {\bf Bound for $U$.} Using Proposition \ref{var22}, we obtain
\begin{eqnarray}\label{chou}
U \le C 2^{\tau}\frac{1}{n}\le C\left( \frac{\ln n}{n}\right)^{2s/(2s+1)} .
\end{eqnarray}

\item {\bf Bound for $W$.} For $q\ge 1$ and $p \ge 2$, we have $g_{\ell}\in \mathbf{B}^s_{p,q}(M)\subseteq \mathbf{B}^s_{2,\infty}(M)$. Hence, by \eqref{berti},
\begin{eqnarray*}
W   \le   C\sum_{j=j_1+1}^{\infty}2^{-2js}\le C2^{-2j_1s}\le C \left(\frac{(\ln n)^3}{n}\right)^{2s}\le   C \left( \frac{\ln n}{n}\right)^{2s/(2s+1)} .
\end{eqnarray*}
For $q\ge 1$ and  $p \in [1,2)$, we have $g_{\ell}\in {B}^{s}_{p,q} (M)\subseteq
{B}^{s+1/2-1/p}_{2,\infty} (M)$. Since $s>1/p$, we have $s+1/2-{1}/{p}>s/(2s+1)$. So, by \eqref{berti},
 \begin{eqnarray*}
W& \le & C\sum_{j= j_1+1}^{\infty} 2^{-2j(s+1/2-{1}/{p})}\le  C 2^{-2j_1(s+1/2-{1}/{p})}\nonumber \\
& \le
&  C\left( \frac{(\ln n)^3}{n}\right)^{2(s+1/2-{1}/{p})}\le C\left( \frac{\ln n}{n}\right)^{2s/(2s+1)}.
\end{eqnarray*}
Hence, for $q\ge 1$, \{$p \ge 2$ and $s>0$\} or \{$p \in [1,2)$ and $s>1/p$\}, we have
 \begin{eqnarray}\label{amp}
W\le   C \left( \frac{\ln n}{n}\right)^{2s/(2s+1)}.
\end{eqnarray}

\item {\bf Bound for $V$.} We have
\begin{eqnarray}\label{dec}
V =V_{1}+V_{2}+V_{3} +V_{4},
\end{eqnarray}
where
$$V_1=\sum_{j=\tau}^{j_1}\sum_{k=0}^{2^j-1} \mathbb{E}\left( (\widehat b_{j,k,\ell}-b_{j,k,\ell})^2{\bf 1}_{\left \lbrace | \widehat b_{j,k,\ell}|\ge\kappa\lambda_n\right \rbrace } {\bf 1}_{\left \lbrace | b_{j,k,\ell}|<\kappa\lambda_n/2\right \rbrace }\right),$$
$$V_2=\sum_{j=\tau}^{j_1}\sum_{k=0}^{2^j-1}  \mathbb{E}\left((\widehat b_{j,k,\ell}-b_{j,k,\ell})^2{\bf 1}_{\left \lbrace | \widehat b_{j,k,\ell}|\ge\kappa\lambda_n\right \rbrace } {\bf 1}_{\left \lbrace |b_{j,k,\ell}|\ge \kappa\lambda_n/2\right \rbrace }\right),$$
$$V_{3}=\sum_{j=\tau}^{j_1}\sum_{k=0}^{2^j-1}\mathbb{E}\left(b_{j,k,\ell}^2{\bf 1}_{\left \lbrace | \widehat b_{j,k,\ell}|<\kappa\lambda_n\right \rbrace } {\bf 1}_{\left \lbrace | b_{j,k,\ell}|\ge 2{\kappa\lambda_n}\right \rbrace }\right)$$
and
$$V_{4}=\sum_{j=\tau}^{j_1}\sum_{k=0}^{2^j-1}\mathbb{E}\left(b_{j,k,\ell}^2{\bf 1}_{\left \lbrace |\widehat  b_{j,k,\ell}|<\kappa\lambda_n\right \rbrace } {\bf 1}_{\left \lbrace | b_{j,k,\ell}|<2{\kappa\lambda_n}\right \rbrace }\right).$$

\begin{itemize}
\item {\bf Bounds for $V_1$ and $V_3$.} The following inclusions hold:

$\left\lbrace |\widehat b_{j,k,\ell}|<\kappa\lambda_n, \ |b_{j,k,\ell}|\ge 2\kappa\lambda_n\right\rbrace   \subseteq \left\lbrace |\widehat b_{j,k,\ell}-b_{j,k,\ell}|>\kappa\lambda_n/2 \right\rbrace$,

$\left\lbrace |\widehat b_{j,k,\ell}|\ge \kappa\lambda_n, \ |b_{j,k,\ell}|<\kappa\lambda_n/2 \right\rbrace \subseteq \left\lbrace |\widehat b_{j,k,\ell}-b_{j,k,\ell}|>\kappa\lambda_n/2 \right\rbrace$

and $\left\lbrace |\widehat b_{j,k,\ell}|<\kappa\lambda_n, \ |b_{j,k,\ell}|\ge 2\kappa\lambda_n \right\rbrace \subseteq \left\lbrace |b_{j,k,\ell}|\le 2|\widehat b_{j,k,\ell}-b_{j,k,\ell}| \right\rbrace$.

So
$$\max(V_1,V_3)\le C \sum_{j=\tau}^{j_1}\sum_{k=0}^{2^j-1}\mathbb{E}\left((\widehat b_{j,k,\ell}-b_{j,k,\ell})^2{\bf 1}_{\left \lbrace |\widehat b_{j,k,\ell}-b_{j,k,\ell}|>\kappa\lambda_n/2\right \rbrace}  \right) .$$
Applying the Cauchy-Schwarz inequality and using Propositions \ref{var2}, \ref{var3} and $2^j\le n$, we have
\begin{eqnarray*}
\lefteqn{\mathbb{E}\left((\widehat b_{j,k,\ell}-b_{j,k,\ell})^2{\bf 1}_{\left \lbrace |\widehat b_{j,k,\ell}-b_{j,k,\ell}|>\kappa\lambda_n/2\right \rbrace}  \right) } & & \nonumber \\
& \le
 & \left(\mathbb{E}\left((\widehat b_{j,k,\ell}-b_{j,k,\ell})^{4}\right)\right)^{1/2}\left(\mathbb{P}\left( |\widehat b_{j,k,\ell}-b_{j,k,\ell}|>\kappa\lambda_n/2\right)\right)^{{1}/{2}} \nonumber \\
& \le
& C\left(\frac{2^j}{n} \right)^{1/2}\left(\frac{1}{n^4} \right)^{1/2}\le C\frac{1}{n^2}.
\end{eqnarray*}
Therefore
\begin{eqnarray}\label{ko}
\max(V_1,V_3) \le    C \frac{1}{n^2}\sum_{j=\tau}^{j_1} 2^{j} \le   C \frac{1}{n^2} 2^{j_1}   \le  C \frac{1}{n} \le C  \left( \frac{\ln n}{n}\right)^{2s/(2s+1)}.
\end{eqnarray}

\item \noindent {\bf Bound for $V_2$.} Using Proposition \ref{var22}, we obtain
\begin{eqnarray*}
\mathbb{E}\left( (\widehat b_{j,k,\ell}-b_{j,k,\ell})^2\right) \le  C\frac{1}{n}\le C\frac{\ln n}{n}.
\end{eqnarray*}
Hence
\begin{eqnarray*}
V_2 \le  C\frac{\ln n}{n}\sum_{j=\tau}^{j_1}\sum_{k=0}^{2^j-1} {\bf 1}_{\left \lbrace |b_{j,k,\ell}|>\kappa\lambda_n/2\right \rbrace}.
\end{eqnarray*}
Let $j_2$ be the integer defined by
\begin{eqnarray}\label{mopp}
 2^{j_2}=\left\lbrack\left(\frac{n}{\ln n}\right)^{1/(2s+1)}\right\rbrack.
\end{eqnarray}
 We have
\begin{eqnarray*}
V_2 \le  V_{2,1}+V_{2,2},
\end{eqnarray*}
where
$$V_{2,1}=C\frac{\ln n}{n}\sum_{j=\tau}^{j_2}\sum_{k=0}^{2^j-1} {\bf 1}_{\left \lbrace |b_{j,k,\ell}|>\kappa\lambda_n/2\right \rbrace}$$
and
$$V_{2,2}=C \frac{\ln n}{n}\sum_{j=j_2+1}^{j_1}\sum_{k=0}^{2^j-1} {\bf 1}_{\left \lbrace |b_{j,k,\ell}|>\kappa\lambda_n/2\right \rbrace}.$$
We have
\begin{eqnarray*}
V_{2,1} \le  C \frac{\ln n}{n}\sum_{j=\tau}^{j_2}2^{j}\le C \frac{\ln n}{n}2^{j_2 }\le C \left(\frac{\ln n}{n}\right)^{2s/(2s+1)}.
\end{eqnarray*}
For $q\ge 1$  and $p\ge 2$, we have $g_{\ell}\in B_{p,q}^s(M)\subseteq \mathbf{B}^s_{2,\infty}(M)$. So, by \eqref{berti},
\begin{eqnarray*}
V_{2,2} & \le  &  C \frac{\ln n}{n \lambda_n^2}\sum_{j=j_2+1}^{j_1}\sum_{k=0}^{2^j-1} b_{j,k,\ell}^2\le C\sum_{j=j_2+1}^{\infty}\sum_{k=0}^{2^j-1} \beta_{j,k}^2\le C 2^{-2j_2s}\\
& \le
& C\left(\frac{\ln n}{n}\right)^{2s/(2s+1)}.
\end{eqnarray*}
For $q\ge 1$, $p \in [1,2)$ and $s>1/p$, using \eqref{berti}, ${\bf 1}_{\left \lbrace |b_{j,k,\ell}|>\kappa\lambda_n/2\right \rbrace}\le C |b_{j,k,\ell}|^{p}/\lambda_n^p=C |\beta_{j,k}|^{p}/\lambda_n^p$ and $(2s+1)(2-p)/2+(s+1/2-1/p) p=2s$, we have
\begin{eqnarray*}
V_{2,2} &  \le & C \frac{\ln n}{n \lambda_n^p}\sum_{j=j_2+1}^{j_1}\sum_{k=0}^{2^j-1}  |\beta_{j,k}|^{p}\le  C \left(\frac{\ln n}{n}\right)^{(2-p)/2}\sum_{j=j_2+1}^{\infty}2^{-j(s+{1}/{2}-{1}/{p})p}\\
&  \le &  C \left(\frac{\ln n}{n}\right)^{(2-p)/2} 2^{-j_2(s+{1}/{2}-{1}/{p})p} \le C \left(\frac{\ln n}{n}\right)^{2s/(2s+1)}.
\end{eqnarray*}
So, for $r\ge 1$, \{$p \ge 2$ and $s>0$\} or \{$p \in [1,2)$ and $s>1/p$\}, we have
\begin{eqnarray}\label{sis}
V_2\le  C  \left(\frac{\ln n}{n}\right)^{2s/(2s+1)}.
\end{eqnarray}

\item \noindent {\bf Bound for $V_4$.} We have
$$V_4 \le \sum_{j=\tau}^{j_1}\sum_{k=0}^{2^j-1} b_{j,k,\ell}^2 {\bf 1}_{\left \lbrace | b_{j,k,\ell}|<2{\kappa\lambda_n}\right \rbrace }.$$
Let $j_2$ be the integer \eqref{mopp}. Then
$$V_4 \le V_{4,1}+V_{4,2},$$
where
$$V_{4,1}= \sum_{j=\tau}^{j_2}\sum_{k=0}^{2^j-1} b_{j,k,\ell}^2 {\bf 1}_{\left \lbrace | b_{j,k,\ell}|<2{\kappa\lambda_n}\right \rbrace }, \ \ \ \ \ \ V_{4,2}= \sum_{j=j_2+1}^{j_1}\sum_{k=0}^{2^j-1} b_{j,k,\ell}^2 {\bf 1}_{\left \lbrace | b_{j,k,\ell}|<2{\kappa\lambda_n}\right \rbrace }.$$
We have
\begin{eqnarray*}
V_{4,1}  \le  C \sum_{j=\tau}^{j_2}2^j \lambda_n^2=C\frac{\ln n}{n} \sum_{j=\tau}^{j_2}2^{j}\le C \frac{\ln n}{n}2^{j_2} \le C \left(\frac{\ln n}{n} \right)^{2s/(2s+1)}.
\end{eqnarray*}
For $q\ge 1$ and $p\ge 2$, we have $g_{\ell}\in B_{p,q}^s(M)\subseteq \mathbf{B}^s_{2,\infty}(M)$. Hence, by \eqref{berti},
$$V_{4,2}\le  \sum_{j=j_2+1}^{\infty}\sum_{k=0}^{2^j-1} \beta_{j,k}^2\le C2^{-2j_2s}\le C \left(\frac{\ln n}{n} \right)^{2s/(2s+1)}.$$
For $q\ge 1$, $p \in [1,2)$ and $s>1/p$, using \eqref{berti}, $b_{j,k,\ell}^2 {\bf 1}_{\left \lbrace | b_{j,k,\ell}|<2{\kappa\lambda_n}\right \rbrace }\le C \lambda_n^{2-p}  |b_{j,k,\ell}|^{p}=C \lambda_n^{2-p}  |\beta_{j,k}|^{p}$ and $(2s+1)(2-p)/2+(s+1/2-1/p) p=2s$, we have
\begin{eqnarray*}
V_{4,2}& \le
& C \lambda_n^{2-p} \sum_{j=j_2+1}^{j_1}\sum_{k=0}^{2^j-1} |\beta_{j,k}|^{p} =C\left(\frac{\ln n}{n}\right)^{(2-p)/2}  \sum_{j=j_2+1}^{j_1}  \sum_{k=0}^{2^j-1} |\beta_{j,k}|^{p} \\
& \le
& C\left(\frac{\ln n}{n}\right)^{(2-p)/2} \sum_{j=j_2+1}^{\infty}2^{-j(s+{1}/{2}-{1}/{p})p}\\
& \le
&  C \left(\frac{\ln n}{n}\right)^{(2-p)/2}2^{-j_2(s+{1}/{2}-{1}/{p})p} \le C\left(\frac{\ln n}{n}\right)^{2s/(2s+1)}.
\end{eqnarray*}
Thus, for $q\ge 1$, \{$p \ge 2$ and $s>0$\} or \{$p \in [1,2)$ and $s>1/p$\}, we have
\begin{eqnarray}\label{qua}
V_4\le C \left(\frac{\ln n}{n}\right)^{2s/(2s+1)}.
\end{eqnarray}
\end{itemize}

It follows from \eqref{dec}, \eqref{ko}, \eqref{sis} and \eqref{qua} that
\begin{eqnarray}\label{benj}
V\le C \left(\frac{\ln n}{n}\right)^{2s/(2s+1)}.
\end{eqnarray}

\end{enumerate}

Combining \eqref{rou}, \eqref{del}, \eqref{chou}, \eqref{amp} and \eqref{benj}, we have, for $q\ge 1$, \{$p \ge 2$ and $s>0$\} or \{$p \in [1,2)$ and $s>1/p$\},
$$\mathbb{E}\left(\int_{0}^{1}(\widehat g_{\ell}(x)-g_{\ell}(x))^2dx\right)\le C \left(\frac{\ln n}{n}\right)^{2s/(2s+1)}.$$
The proof of Theorem \ref{theo2} is complete.

\endproof

\paragraph{Acknowledgement.}
This work is supported by ANR grant NatImages, ANR-08-EMER-009.

\end{document}